\documentclass[12pt]{amsart}
\usepackage{times}
\usepackage{amsfonts}
\usepackage{amsthm}
\usepackage{amsmath}
\advance \topmargin by -\headheight
\advance \topmargin by -\headsep
\evensidemargin \oddsidemargin
\newtheorem{theorem}{Theorem}[section]
\newtheorem{lemma}[theorem]{Lemma}

\newtheorem{remark}[theorem]{Remark}

\newtheorem{definition}[theorem]{Definition}

\begin{document}

\title{Amenability, Tubularity, and Embeddings into $\mathcal R^{\omega}$}

\author{Kenley Jung}

\dedicatory{Dedicated to Ed Effros on the occasion of his 70th birthday}

\address{Department of Mathematics, University of California,
Los Angeles, CA 90095-1555,USA}

\email{kjung@math.ucla.edu}
\thanks{Research supported in part by the NSF}

\begin{abstract} Suppose $M$ is a tracial von Neumann algebra embeddable 
into $\mathcal R^{\omega}$ (the ultraproduct of the hyperfinite 
$\mathrm{II}_1$-factor) and $X$ is an $n$-tuple of selfadjoint generators 
for $M$.  
Denote by $\Gamma(X;m,k,\gamma)$ the microstate space of $X$ of order 
$(m,k,\gamma)$.  We say that $X$ is \textit{tubular} if for any $\epsilon 
>0$ there exist $m \in \mathbb N$ and $\gamma>0$ such that if 
$(x_1,\ldots, x_n), (y_1, \ldots, y_n) \in \Gamma(X;m,k,\gamma),$ then 
there exists a $k \times k$ unitary $u$ satisfying $|ux_iu^* - y_i|_2 < 
\epsilon$ for each $1 \leq i \leq n.$ We show that the following 
conditions are equivalent:

\vspace{.1in}

\begin{itemize}
\item $M$ is amenable (i.e., injective).
\item $X$ is tubular.
\item Any two embeddings of $M$ into $\mathcal R^{\omega}$ are conjugate by a unitary $u \in \mathcal R^{\omega}$.
\end{itemize}
\end{abstract}
\maketitle

\section{Introduction}

One version of Voiculescu's free entropy involves \textit{matricial 
microstates}.  Given an $n$-tuple of selfadjoint elements $X = 
\{x_1,\ldots, x_n\}$ in a von Neumann algebra $M$ with trace $\varphi$, 
denote by $\Gamma(X;m,k,\gamma)$ the set of all $n$-tuples of selfadjoint 
$k \times k$ matrices $(a_1,\ldots,a_n)$ such that for any $1 \leq p \leq 
m$ and $1\leq i_1, \ldots, i_p \leq n,$

\[ |tr_k(a_{i_1} \cdots a_{i_p}) - \varphi(x_{i_1} \cdots x_{i_p})| < \gamma
\]

\noindent where $tr_k$ denotes the normalized trace on the $k \times k$ 
matrices.  Such an $n$-tuple $(a_1, \ldots, a_n)$ is called a microstate 
for $X$ and the sets $\Gamma(X;m,k,\gamma)$ are called the microstate 
spaces of $X$.

One can use microstates and random matrices to obtain nonisomorphism 
results (for an overview see [9]).  At some point all of these kinds of 
arguments rely on the fact that amenability 
(=hyperfiniteness=injectivity by [1]; see also [2]-[6] and [8]) 
forces an extremely rigid condition 
on the microstate spaces.  When $X^{\prime \prime}$ is amenable then the 
following is true: given $\epsilon >0$ there exist $m \in \mathbb N$ and 
$\gamma >0$ such that for any two elements $\xi = (\xi_1, \ldots, \xi_n), 
\eta = (\eta_1, \ldots, \eta_n) \in \Gamma(X;m,k,\gamma)$ there exists a 
$k \times k$ unitary $u$ such that for each $1 \leq i \leq n,$ $|u\xi_i 
u^* - \eta_i|_2 < \epsilon.$ Roughly speaking then, a microstate space of 
$X$ is the neighborhood of the unitary orbit of a single microstate of the 
space.

In this note we observe the converse.  To be more exact, suppose 
$X^{\prime \prime}$ embeds into $\mathcal R^{\omega}$ and satisfies the 
aforementioned geometric property: for any $\epsilon >0$ there exist $m 
\in \mathbb N$ and $\gamma >0$ such that if $\xi = (\xi_1, \ldots, \xi_n), 
\eta = ( \eta_1, \ldots, \eta_n) \in \Gamma(X;m,k,\gamma),$ then there 
exists a $k \times k$ unitary $u$ satisfying $|u \xi_i u^* - \eta_i|_2 < 
\epsilon, 1 \leq i \leq n.$ We show that $X^{\prime \prime}$ is 
semidiscrete in the sense of [5] and thus, by [2] amenable.  

We also notice that this characterization is equivalent to the condition 
that any two embeddings $\sigma, \pi$ of $X^{\prime \prime}$ into 
$\mathcal R^{\omega}$ are conjugate by a unitary $u \in \mathcal 
R^{\omega}$.  In the course of the proof we will introduce seemingly 
weaker notions of tubularity: finite tubularity, quasitubularity, and 
finite quasitubularity.  All of these notions will be shown to coincide 
with amenability as well. 

\section{Some other characterizations of Amenability}

Throughout $M$ will be a von Neumann algebra with tracial, faithful, 
normal state $\varphi$ and $X = \{x_1,\ldots, x_n\}$ is a finite 
generating set of selfadjoints of $M$.  Fix a nontrivial ultrafilter 
$\omega$ of $\mathbb N$; there exists an obvious net $i : \Lambda 
\rightarrow \mathbb N$ such that for any bounded sequence of complex 
numbers, $\langle c_n \rangle_{n=1}^{\infty}$, $\omega(\langle c_n 
\rangle_{n=1}^{\infty}) = \lim_{\lambda \in \Lambda} c_{i (\lambda) }.$ 
$\mathcal R$ will be a fixed copy of the hyperfinite 
$\mathrm{II}_1$-factor, and $\mathcal R^{\omega}$ will be the associated 
ultraproduct construction with respect to $\omega$ and $\mathcal R$.  $Q : 
\ell^{\infty}(\mathcal R) \rightarrow \mathcal R^{\omega}$ is the obvious 
quotient map.  By "embedding" we mean a normal, injective $*$-homorphism 
which preserves units and the traces.

We will be constantly working with finite tuples of operators and for this reason introduce the following, somewhat abusive notation.  For an ordered $n$-tuple $\xi = \{\xi_1,\ldots, 
\xi_n\}$ in a von Neumann algebra $N$ and a unitary $u \in N$ $u\xi u^* = 
\{ u\xi_1 u^*, \ldots, u \xi_n u^*\}$.  If $F$ is a map from $N$ into some 
set $S$, then we carelessly write $F(\xi)$ for $\{ F(\xi_1), \ldots, F(\xi_n)\}$. If $\eta = 
\{\eta_1, \ldots, \eta_n\} \subset N$, then $\xi - \eta = \{\xi_1 - 
\eta_1, \ldots, \xi_n - \eta_n\}$ and $|\xi|_2 = (\sum_{i=1}^n 
\tau(\xi_i^2))^{\frac{1}{2}}$ where $\tau$ is a tracial state on $N$.  We 
will assume throughout that $M$ embeds into $\mathcal R^{\omega}$.  For $k 
\in \mathbb N,$ $M_k(\mathbb C)$ denotes the $k\times k$ matrices and 
$U_k$ is the unitary group of $M_k(\mathbb C).$

\begin{definition} $X$ is $N$-tubular ($N \in \mathbb N$) if for any 
$\epsilon >0$ there exist $m \in \mathbb N$ and $\gamma >0$, such that if 
$\xi_1, \ldots, \xi_{N+1} \in \Gamma(X;m,k,\gamma)$, then there exists a 
$u \in U_k$ satisfying $|u\xi_i u^* - \xi_j|_2 < \epsilon$ for some $1 
\leq i < j \leq N+1$.  $X$ is finitely tubular if $X$ is $N$-tubular for 
some $N \in \mathbb N$.  $X$ is simply tubular if $X$ is $1$-tubular. 
\end{definition}

Thus, $X$ is finitely tubular if it can be encapsulated in the unitary 
orbits of no more than $N$ of its microstates.  Clearly, tubularity 
coincides with the definition given in the introduction.  We have a 
similar notion for quasitubularity:

\begin{definition} $X$ is $N$-quasitubular ($N \in \mathbb N$) if for any 
$\epsilon >0$ there exist $m \in \mathbb N$ and $\gamma >0$, such that for 
any $\xi_1, \ldots, \xi_{N+1} \in \Gamma(X;m,k,\gamma)$ there exists a $p$ 
(dependent on $\xi_1, \ldots, \xi_{N+1}$) and a unitary $u$ in 
$M_k(\mathbb C) \otimes M_p(\mathbb C)$ satisfying $|u(\xi_i \otimes 
I_{p})u^* - \xi_j \otimes I_p |_2 < \epsilon$ for some $1 \leq i < j \leq 
N+1$.  $X$ is finitely quasitubular if $X$ is $N$-quasitubular for some 
$N$.  $X$ is simply quasitubular if $X$ is $1$-quasitubular.  
\end{definition}

\begin{remark} Obviously if $X$ is $N$-tubular, then $X$ is 
$N$-quasitubular. \end{remark}

\begin{lemma} If $X$ is tubular, then any two embeddings $\sigma, \pi$ of 
$X^{\prime \prime}$ into $\mathcal R^{\omega}$ are conjugate by a unitary 
in $\mathcal R^{\omega}$. \end{lemma}

\begin{proof} Suppose $X$ is tubular and $\sigma, \pi : X^{\prime \prime} 
\rightarrow \mathcal R^{\omega}$ are two embeddings.  Find algebras $A_k 
\subset \mathcal R$ such that for each $k$, $A_k \simeq M_{r(k)}(\mathbb 
C)$ for some $r(k) \in \mathbb N$ and such that for each $1 \leq j \leq n$ 
there exist sequences $y_j = \langle y_{jk} \rangle_{k=1}^{\infty}$, $z_j 
= \langle z_{jk} \rangle_{k=1}^{\infty} \in \ell^{\infty}(\mathcal R)$ 
satisfying $y_{jk}, z_{jk} \in A_k$ for each $k$, $\pi(x_j) = Q (y_j),$ 
and $\sigma(x_j) = Q(z_j).$

   For each $p \in \mathbb N$ there exists by tubularity a corresponding 
$m(p) \in \mathbb N$ such that for any $k \in \mathbb N$, if $\xi, \eta 
\in \Gamma(F;m(p),k, m(p)^{-1}),$ then there exists a $u \in U_k$ 
satisfying $|u \xi u^* - \eta|_2 < p^{-1}.$ For each $k$ pick a unitary 
$w_k \in A_k \simeq M_{r(k)}(\mathbb C)$ satisfying

\[ \max_{1 \leq j \leq n} |w_k y_{jk} w_k^* - 
z_{jk}|_2 = \inf_{w \in U(A_k)} \left( \max_{1 \leq j \leq n} |w y_{jk} 
w^* - z_{jk}|_2 \right) \] 

\noindent where $U(A_k)$ ($\simeq U_k$) is the unitary group of 
$A_k$.  Define $w = \langle w_k \rangle_{k=1}^{\infty}$; $w \in 
\ell^{\infty}(M).$ Given an $N \in \mathbb N$, for $\lambda$ sufficiently 
large $(y_{1i(\lambda)}, \ldots, y_{ni(\lambda)}), (z_{1i(\lambda)}, 
\ldots, z_{ni(\lambda)}) \in \Gamma(F;m(N), i(\lambda), m(N)^{-1})$ which 
in turn implies that for such $\lambda$ and all $1 \leq j \leq n$
   
\[ |w_{i(\lambda)}y_{i(\lambda)j} w_{i(\lambda)}^* - z_{i(\lambda)j}|_2 < 
N^{-1}. \]

\noindent Set $u = Q(w) \in \mathcal R^{\omega}$. We have that for all $1 
\leq j \leq n,$ $\sigma(x_j) = u \pi(x_j) u^*$ which completes the proof. 
\end{proof}

\begin{lemma} If any two embeddings $\sigma$ and $\pi$ of $M$ into 
$\mathcal R^{\omega}$ are conjugate by a unitary in $\mathcal R^{\omega}$, 
then $X$ is quasitubular. \end{lemma}

\begin{proof} Suppose by contradiction that $X$ is not quasitubular.  For 
some $\epsilon_0 >0$ and any $m \in \mathbb N$ and $\gamma >0$ there 
exists an $N \in \mathbb N$ and $\xi, \eta \in \Gamma(X;m,N,\gamma)$ such 
that for all $p \in \mathbb N,$

\[ \inf_{u \in U_{Np}} |u (\xi \otimes I_p) u^* - \eta \otimes I_p|_2 > 
\epsilon_0.\]

\noindent Thus, for each $m \in \mathbb N$ we can find a corresponding 
$N_m \in \mathbb N$ and $\xi_m, \eta_m \in \Gamma(X;m, N_m, m^{-1})$ such 
that for any $k \in \mathbb N$ $\inf_{u \in U_{kN_m}} |u (\xi_m \otimes 
I_k ) u^* - \eta_m \otimes I_k |_2 > \epsilon_0.$ Without loss of 
generality we may assume that the operator norms of any of the elements in 
$\xi_m$ or $\eta_m$ are strictly less than $C = \max_{x \in X} \|x\| +1$.  
For each $m$, $\mathcal R = M_{N_m}(\mathbb C) \otimes \mathcal R_m$ where 
$\mathcal R_m \simeq \mathcal R$; define $x_m = \xi_m \otimes I$ and $y_m 
= \eta_m \otimes I$ with respect to this decomposition of $\mathcal R$ and 
set $x = \langle x_m \rangle_{m=1}^{\infty}$ and $y = \langle y_m 
\rangle_{m=1}^{\infty}$.  It is not too hard to see that we can find two 
embeddings $\pi, \sigma: M \rightarrow \mathcal R^{\omega}$ satisfying 
$\pi(X) = Q(x)$ and $\sigma(X) = Q(y)$.

By hypothesis there exists a unitary $u \in \mathcal R^{\omega}$ 
satisfying $\sigma(x) = u \pi(x) u^*$ for all $x \in M$.  We can find some 
$w = \langle w_m \rangle_{m=1}^{\infty} \in \ell^{\infty}(\mathcal R)$ 
such that $Q(w) = u$.  Because $u$ is a unitary we can assume that for 
each $m$, $w_m \in \mathcal R$ is a unitary. The condition $\sigma(x) = u 
\pi(x) u^*$ implies that there exists a $\lambda \in \Lambda$ such that

\[ |w_{i(\lambda)} (\xi_{i(\lambda)} \otimes I) w_{i(\lambda)}^* - 
\eta_{i(\lambda)} \otimes I |_2 = |w_{i(\lambda)} x_{i(\lambda)} 
w_{i(\lambda)}^* - y_{i(\lambda)}|_2 < \epsilon_0/3C.\]

\noindent Now, $w_{i(\lambda)} \in \mathcal R = M_{N_{i(\lambda)}}(\mathbb 
C) \otimes R_{i(\lambda)}$ and by standard approximations we can find some 
$p_0 \in \mathbb N$, a unital $*$-algebra $A_{\lambda} \simeq 
M_{p_0}(\mathbb C)$ with $A_{\lambda} \subset \mathcal R_{i(\lambda)}$, 
and a unitary $u \in M_{N_{i(\lambda)}} \otimes A_{\lambda} \subset 
\mathcal R$ satisfying $|u - w_{i(\lambda)}|_2 < \epsilon_0/3C.$ It is 
then clear that

\[ \inf_{v \in U_{i(\lambda)p_0}} | v(\xi_{i(\lambda)} \otimes I_{p_0}) v^* - \eta_{i(\lambda)} \otimes I_{p_0}|_2 \leq |u(\xi_{i(\lambda)} \otimes I_{p_0}) u^* - \eta_{i(\lambda)} \otimes I_{p_0}|_2 < \epsilon_0. \] 

\noindent which contradicts the initial hypothesis.
\end{proof}

We now present a lemma which is undoubtedly known but which we will prove 
for completeness.  Recall that a von Neumann algebra $N$ is semidiscrete 
if there exist nets $\langle \phi_{j} \rangle_{j \in \Omega}$, and 
$\langle \psi_j \rangle_{j \in \Omega}$ of unital completely positive maps 
$\phi_j : N \rightarrow M_{n_j}(\mathbb C)$, $\psi_j : M_{n_j}(\mathbb C) 
\rightarrow N$, $n_j \in \mathbb N$ such that for any $x \in N$, $(\phi_j 
\circ \psi_j)(x) \rightarrow x$ $\sigma$-weakly.  This is not the original definition of semidiscreteness found in [5] (which demands that the maps only have finite rank), but it is equivalent to that definition by [3] and [4].  Semidiscreteness, which 
was introduced in [5], is yet another characterization of amenability ([1] 
again).

\begin{lemma} Suppose $A$ is a tracial von Neumann algebra. Assume that 
for some finite set of generators $F$ of $A$ and any $\epsilon >0$ there 
exist an embedding $\pi_{\epsilon}$ of $A$ into a tracial von Neumann 
algebra $A_{\epsilon}$ and a unital, finite dimensional algebra 
$B_{\epsilon} \subset A_{\epsilon}$ with the property that every element 
of $\pi_{\epsilon}(F)$ is contained in the $\epsilon$-neighborhood of 
$B_{\epsilon}$ with respect to the $|\cdot |_2$-norm of $A_{\epsilon}$.  
Then $A$ is semidiscrete, and thus amenable. \end{lemma}

\begin{proof} The hypothesis implies that for any finite set $S$ of $A$ 
and any $\epsilon >0$ there exists an embedding $\pi_{\epsilon}$ of $A$ 
into a tracial von Neumann algebra $A_{\epsilon}$ and a finite dimensional 
unital subalgebra $B_{\epsilon}$ of $A_{\epsilon}$ such that every element 
of $\pi_{\epsilon}(S)$ is contained in the $\epsilon$-neighborhood of 
$B_{\epsilon}$ with respect to $|\cdot|_2$.  This is because the elements 
of $\pi_{\epsilon}(F)$ can be approximated in $|\cdot|_2$-norm by elements 
in $B_{\epsilon}$ with operator norms no bigger than the maximum of the 
operator norms of elements in $F$ (one uses conditional expectations as we 
do below) and because multiplication is $|\cdot|_2$-continuous on operator 
norm bounded sets.  This gives the implication for $S$ consisting of 
polynomials of elements from $F$ and the general case follows immediately.

To complete the proof it suffices to show that for any finite $S \subset 
A$ we can construct sequences of u.c.p. maps $\phi_n : A \rightarrow 
M_{k_n}(\mathbb C)$ and $\psi_n: M_{k_n}(\mathbb C) \rightarrow A$, $k_n 
\in \mathbb N$, such that for any normal linear functional $f$ on $N,$ 
$f((\phi_n \circ \psi_n)(x)) \rightarrow f(x)$ as $n \rightarrow \infty$.  
Thus, let the finite subset $S$ of $A$ be given.  By the first paragraph 
for each $n \in \mathbb N$ we can find a tracial von Neumann algebra 
$A_n$, an embedding $\pi_n : A \rightarrow A_n$ and a finite dimensional, 
unital subalgebra $B_n$ of $A_n$ such that every element of $S$ is 
contained in the $n^{-1}$-neighborhood of $B_n$ with respect to the 
$|\cdot|_2$-norm of $A_n$.  Define $\phi_n : A \rightarrow B_n$ to be the 
composition of $\pi_n$ with the conditional expectation $E_n$ of $A_n$ 
onto $B_n$. Define $F_n$ to be the conditional expection of $A_n$ onto 
$\pi_n(A)$ restricted to $B_n$; so $F_n : B_n \rightarrow \pi_n(A)$.  
Define $\psi_n : B_n \rightarrow A$ to be the composition of $F_n$ with 
$\pi_n^{-1}$.  Obviously for any $x \in S$

\[ |(\psi_n \circ \phi_n)(x) - x|_2 = |F_n(E_n(\pi_n(x))) - \pi_n(x)|_2 \leq |E_n(\pi_n(x)) - \pi_n(x)|_2 <  n^{-1}.\]

\noindent Now $\langle (\psi_n \circ \phi_n)(x) \rangle_{n=1}^{\infty}$ is 
a sequence in $A$ uniformly bounded in the operator norm and $(\psi_n 
\circ \phi_n)(x) \rightarrow x$ in the $|\cdot|_2$-norm.  This implies 
$(\psi_n \circ \phi_n)(x) \rightarrow x$ $\sigma$-weakly for every $x \in 
S.$

\end{proof} 

\begin{remark} Suppose $N$ is a von Neumann algebra, and $A \subset N$ is a von Neumann subalgebra. By the above lemma if for some finite set of generators $F$ of $A$ and any $\epsilon >0$ there exists a finite dimensional algebra $B \subset N$ such that every element of $F$ is contained in the $\epsilon$-neighborhood of $B$ with respect to $|\cdot |_2$, then $A$ is amenable.
\end{remark}

\begin{lemma} If $X$ is quasitubular, then $M$ is amenable.
\end{lemma}

\begin{proof} For each $m \in \mathbb N$ there exist a $k(m) \in \mathbb N$ and a $\xi_m = (\xi_{1m}, \ldots, \xi_{nm}) \in \Gamma(X;m,k(m), m^{-1})$.  $\mathcal R = \mathcal R_1 \otimes \mathcal R_2 \otimes \mathcal R_3$ where $\mathcal R_i \simeq \mathcal R$, $i=1, 2, 3.$  Define for $1\leq j \leq n$, $y_j = \langle I \otimes \xi_{jm} \otimes I \rangle_{m=1}^{\infty} \in \ell^{\infty}(\mathcal R)$ where $I \otimes \xi_{jm} \otimes I \in \mathcal R_1 \otimes \mathcal B_m  \otimes \mathcal R_3 \subset \mathcal R_1 \otimes \mathcal R_2  \otimes \mathcal R_3 = \mathcal R$ and $M_{k(m)}(\mathbb C) \simeq B_m \subset R_2.$  It is clear that there exists a (trace preserving) embedding $\pi : X^{\prime \prime} \rightarrow R^{\omega}$ satisfying $\pi(x_j) = Q(y_j)$.  Given $\epsilon >0$ we can find by quasitubularity an $m_0 \in \mathbb N$ such that for any $k \in \mathbb N$ and $\xi, \eta \in \Gamma(X;m_0, k, m_0^{-1})$ there exists a corresponding $p(k) \in \mathbb N$ and unitary $u$ of $M_k(\mathbb C) \otimes M_{p(k)}(\mathbb C)$ satisfying $|u(\xi \otimes I_{p(k)})u^* - \eta \otimes I_{p(k)}|_2 < \epsilon$.  Fix once and for all, an $N \in \mathbb N$ such that there exists a $\zeta \in \Gamma(X;m_0, N, m_0^{-1}).$  We can regard $\zeta$ as a subset of $\mathcal R_1$ by finding a copy $A$ of the $N \times N$ matrices in $\mathcal R_1$ and for each $m$ consider $\zeta \otimes I \otimes I, I \otimes \xi \otimes I \subset ( A \otimes \mathcal B_m ) \otimes \mathcal R_3 \subset \mathcal R$.  For each $m$ we can find an algebra $D_m$ isomorphic to a full matrix algebra, and a unitary $u_m$ of $A \otimes \mathcal B_m \otimes D_m$ such that

\[ |u_m(\zeta \otimes I \otimes I) u_m^* - I \otimes \xi \otimes I|_2 < \epsilon.
\]

\noindent Define $\rho: A \rightarrow \mathcal R^{\omega}$ by $\rho(x) = Q \langle u_m(x \otimes I \otimes I) u_m^* \rangle_{m=1}^{\infty}$.  We have now shown that for $\epsilon >0$, every element of $\pi(X)$ is within the $| \cdot |_2$ $\epsilon$-ball of the finite dimensional full matrix algebra $\rho(A)$.

By the remark this implies that the von Neumann algebra generated by $\pi(X)$ is semidiscrete, and hence, amenable.  Since $\pi$ is an isomorphism, $X^{\prime \prime} = M$ is amenable.        
\end{proof}

 At this point we have already demonstrated the equivalence claimed in 
the introduction.  We will now prove $\epsilon$ more.  Notice that the 
content of the observation below is the implication that finite 
quasitubularity implies amenability.  This will be very much like the 
proof above modulo some technicalities.  We could have gone straight to 
the following more technical argument without proving Lemma 2.8, but the 
proof of Lemma 2.8 has the advantage of being more lucid.

\begin{lemma} The following are equivalent:

\begin{list}{}{}
\item (1) $M$ is amenable.
\item (2) $X$ is tubular.
\item (3) $X$ is finitely tubular.
\item (4) $X$ is quasitubular.
\item (5) $X$ is finitely quasitubular.
\item (6) Any two embeddings of $M$ into $\mathcal R^{\omega}$ are conjugate by a unitary in $\mathcal R^{\omega}$.

\end{list}
\end{lemma}

\begin{proof} By [2] $M$ is amenable iff $M$ is hyperfinite and thus by 
[7], $X$ must be tubular.  Hence, (1) $\Rightarrow$ (2). By Lemma 2.4, (2) 
$\Rightarrow$ (6) and by Lemma 2.5 (6) $\Rightarrow$ (4).  Clearly (4) 
$\Rightarrow$ (5).

I will now show (5) $\Rightarrow$ (1). Suppose $X$ is finitely 
quasitubular.  Find the smallest $N \in \mathbb N$ for which $X$ is 
$N$-quasitubular.  By Lemma 2.6 we can assume $N >1$. $X$ is not 
$(N-1)$-quasitubular, which implies that we can find some $\epsilon_1 >0$, 
such that for any $m \in \mathbb N$ there exists some $k \in \mathbb N$, 
and $\xi_1, \ldots, \xi_N \in \Gamma(X;m,k, m^{-1})$ such that for any $p 
\in \mathbb N$ and unitary $u$ of $M_{k}(\mathbb C) \otimes M_p(\mathbb 
C)$,

\[ \min_{1 \leq i < j \leq N} |u(\xi_i \otimes I_p)u^* - \xi_j \otimes I_p|_2 > \epsilon_1.
\]

\noindent Let $\epsilon_1 > \epsilon >0$ be given.  There exists an $m_0 
\in \mathbb N$ such that if $\xi_1, \ldots, \xi_{N+1} \in \Gamma(X;m_0, k, 
m_0^{-1})$ then there exists some $p \in \mathbb N$ and unitary $u$ in 
$M_{k_0}(\mathbb C) \otimes M_p(\mathbb C)$ ($u$ and $p$ dependent on the 
$\xi_i$) satisfying $|u(\xi_i \otimes I_p)u^* - \eta|_2 < \epsilon$ for 
some $1 \leq i < j \leq N+1$.  On the other hand, by our initial remark, 
there is some $k(0) \in \mathbb N$, and $\eta_1, \ldots, \eta_N \in 
\Gamma(X;m_0,k(0), m_0^{-1})$ such that for any $p \in \mathbb N$ and 
unitary $u$ of $M_{k(0)}(\mathbb C) \otimes M_p(\mathbb C)$,

\begin{eqnarray} \min_{1 \leq i < j \leq N} |u(\xi_i \otimes I_p)u^* - \xi_j \otimes I_p|_2 > \epsilon_1 > \epsilon.
\end{eqnarray}

Now there exists a sequence $\langle \zeta_m \rangle^{\infty}_{m=1}$ such 
that for each $m$ there exists a $k(m) \in \mathbb N$ with $\zeta_m \in 
\Gamma(X;m,k(m), m^{-1})$.  Identifying $M_{k(0)k(m)}(\mathbb C)$ with 
$M_{k(0)}(\mathbb C) \otimes M_{k(m)}(\mathbb C)$ for all $m \geq m_0$ we 
have that $\xi_1 \otimes I_{k(m)}, \ldots, \xi_N \otimes I_{k(m)}, 
I_{k(0)} \otimes \zeta_m \in \Gamma(X;m_0, k(0)k(m), m_0^{-1})$.  (1) and the 
$N$-quasitubularity of $X$ implies that there must exist an $1 \leq i_m 
\leq N$, a $p_m \in \mathbb N$, and a unitary $u_m$ of $M_{k(0) 
k(m)}(\mathbb C) \otimes M_{p_m}(\mathbb C)$ satisfying

\begin{eqnarray} |u_m(\xi_{i_m} \otimes I_{p_m} )u_m^* - \zeta_m \otimes I_{p_m}|_2 < \epsilon.
\end{eqnarray}

\noindent $\langle i_m \rangle_{m=1}^{\infty}$ is a sequence taking integral values between $1$ and $N$ and thus $i_m =i$ for some $1 \leq i \leq N$ and infinitely many $m$.  Without loss of generality assume $i_{m_q} =1$ for some increasing sequence $\langle m_q \rangle_{q=1}^{\infty}$ of $\mathbb N$.  

For each $q \in \mathbb N$ set $\theta_q = \zeta_{m_q} \in \Gamma(X;m_q,k(m_q), m_q^{-1})$.  $\mathcal R = \mathcal R_1 \otimes \mathcal R_2 \otimes \mathcal R_3$ where $\mathcal R_i \simeq \mathcal R$, $i=1, 2, 3.$  Define $Y = \langle I \otimes \theta_q \otimes I \rangle_{q=1}^{\infty} \in (\ell^{\infty}(\mathcal R))^n$ where $I \otimes \theta_q \otimes I \in \mathcal R_1 \otimes \mathcal B_q  \otimes \mathcal R_3 \subset \mathcal R_1 \otimes \mathcal R_2  \otimes \mathcal R_3 = \mathcal R$ and $M_{k(m_q)}(\mathbb C) \simeq B_q \subset R_2.$  Because $\lim_{q \rightarrow \infty} m_q = \infty$, it is clear that there exists a (trace preserving) embedding $\pi : X^{\prime \prime} \rightarrow R^{\omega}$ satisfying $\pi(X) = Q(Y)$.    We can regard $\xi_1$ as a subset of $\mathcal R_1$ by fixing a copy $A$ of the $k(0) \times k(0)$ matrices in $\mathcal R_1$ and for each $q$ consider $\xi_1 \otimes I \otimes I, I \otimes \theta_q \otimes I \subset ( A \otimes \mathcal B_m ) \otimes \mathcal R_3 \subset \mathcal R$.  For each $q$ (2) provides an algebra $D_q$ isomorphic to a full matrix algebra, and a unitary $v_q$ of $A \otimes \mathcal B_q \otimes D_q$ such that

\[ |v_q(\xi_1 \otimes I \otimes I) v_q^* - I \otimes \theta_q \otimes I|_2 < \epsilon.
\]

\noindent Define $\rho: A \rightarrow \mathcal R^{\omega}$ by $\rho(x) = Q 
\langle v_q(x \otimes I \otimes I) v_q^* \rangle_{q=1}^{\infty}$.  We have 
shown that every element of $\pi(X)$ is within the $| \cdot |_2$ 
$\epsilon$-ball of the finite dimensional full matrix algebra $\rho (A)$.  
Thus, for every $\epsilon >0$ there exists an isomorphic copy of $M$ in 
$\mathcal R^{\omega}$ such that $X$ identified in 
$\mathcal R^{\omega}$ with respect to this embedding is within the 
$|\cdot|_2$-$\epsilon$ ball of a 
finite 
dimensional subalgebra of $\mathcal R^{\omega}$.  It follows that $M$ is 
amenable 
by Lemma 2.6.  (1) follows.

We now have the equivalence of conditions (1), (2), (4), (5) and (6).  To 
conclude, (3) $\Rightarrow$ (5) = (2) $\Rightarrow$ (3) so 
condition (3) is equivalent to all the other conditions as well.
\end{proof}

\begin{remark} It should be fairly obvious to the reader by now that 
conditions (1)-(6) of Lemma 2.9 are equivalent to the condition that there 
exist finitely many embeddings $\pi_1, \ldots, \pi_n$ of $M$ into 
$\mathcal 
R^{\omega}$ such that for any other embedding $\sigma$ of $M$ into 
$\mathcal R^{\omega}$, there exists a $1 \leq i \leq n$ and a unitary $u$ 
in $\mathcal R^{\omega}$ such that $\sigma$ is conjugate to $\pi_i$ via 
$u$. \end{remark}

\noindent{\it Acknowledgements.} I would like to thank Nate Brown, Ed Effros, Sorin 
Popa, and Dimitri Shlyakhtenko for useful conversations.  I am especially grateful to Nate for simplifying my original argument and to Ed for the discussions on injectivity.

\end{document}